\documentclass[11pt]{article}
\usepackage[T1]{fontenc}
\usepackage{lmodern,amsmath,amssymb,amsthm,mathtools}
\usepackage[margin=1in]{geometry}
\usepackage{microtype,enumitem,booktabs}
\usepackage[colorlinks=true,linkcolor=blue!50!black,urlcolor=blue!50!black,citecolor=blue!50!black]{hyperref}
\usepackage{xcolor}
\newtheorem{theorem}{Theorem}[section]
\newtheorem*{externaltheorem}{Theorem}
\newtheorem{proposition}[theorem]{Proposition}
\newtheorem{lemma}[theorem]{Lemma}

\theoremstyle{remark}
\theoremstyle{definition}
\newtheorem{definition}[theorem]{Definition}

\newcommand{\R}{\mathbb R}
\newcommand{\F}{\mathbf F}
\newcommand{\dimH}{\dim_{\mathrm H}}
\newcommand{\supp}{\operatorname{supp}}
\newcommand{\dist}{\operatorname{dist}}
\newcommand{\N}{\mathcal N}
\newcommand{\HH}{\mathcal H}
\newcommand{\dd}{\,\mathrm d}
\newcommand{\eps}{\varepsilon}

\setlist{nosep}
\makeatletter
\renewcommand{\@maketitle}{%
  \begingroup
  \centering
  {\fontsize{13}{17}\usefont{OT1}{cmr}{bx}{n}\MakeUppercase{\@title}\par}
  \vspace{2.0mm}
  {\fontsize{9}{11}\usefont{OT1}{cmr}{m}{n}\MakeUppercase{\@author}\par}
  \endgroup
  \vspace{.8em}}

\newcommand{\introsubsection}{%
  \@startsection{subsection}{2}{\z@}%
    {-3.25ex \@plus -1ex \@minus -.2ex}%
    {-1em}%
    {\normalfont\large\bfseries}}
\newcommand{\proofsubsection}{%
  \@startsection{subsection}{2}{\z@}%
    {-2ex \@plus -.5ex \@minus -.2ex}%
    {.75ex \@plus .1ex}%
    {\normalfont\normalsize}}
\newcommand{\proofstep}{%
  \@startsection{subsubsection}{3}{\z@}%
    {-1.5ex \@plus -.3ex \@minus -.2ex}%
    {1.5ex \@plus .2ex}%
    {\normalfont\normalsize}}
\makeatother
\numberwithin{equation}{section}
\title{Hyperplane Furstenberg Sets in Higher Dimensions}
\author{Peiqi Nie}
\date{}
\usepackage{sectsty}
\sectionfont{\fontsize{14}{16}\selectfont}
\begin{document}
\maketitle
\begin{abstract}
The hyperplane Furstenberg problem asks how small the Hausdorff dimension of \(E\subset\mathbb R^n\) can be when a \(t\)-dimensional family of affine $(n-1)$-dimensional hyperplanes meets \(E\) in sets of dimension at least \(s\). For \(0<s\le n-1\) and \(0<t\le n\), we prove the sharp bound in $\R^n$
\[
\dim_{\mathrm H}E\ge \mathbf F(s,t;n,n-1),
\]
where \(\mathbf F\) is Gan’s Furstenberg index proposed in prime fields. We also construct examples attaining this bound for all parameters. 
Our main reduction theorem, Theorem~3.1, states that for $n\ge3$,
$0<s\le n-2$, and $1<t\le n$, every weighted $(s,t)$ configuration
$(E,H)$ at scale $\delta$ produces a dual
weighted $(s,t-1)$ configuration $(P,\mathcal T)$ inside a
hyperplane $H_{y_0}\cong\mathbb R^{n-1}$. The reduction uses radial projection and energy
estimates for Frostman measures. Iterating this reduction and
applying the Cs\"ornyei--Stull incidence theorem yields the bound.
\end{abstract}

\setcounter{tocdepth}{1}
\begingroup
\large
\makeatletter
\renewcommand{\l@section}[2]{%
  \addvspace{2pt}%
  \@dottedtocline{1}{0em}{1.7em}{#1}{#2}}
\renewcommand{\tableofcontents}{%
  \par\medskip
  {\centering\fontsize{10}{12}\usefont{OT1}{cmr}{m}{sc}\contentsname\par}
  \smallskip
  \@starttoc{toc}}
\tableofcontents
\makeatother
\endgroup

\section{Introduction}\label{sec:introduction}

\introsubsection{Main results.} We prove the following sharp lower bound for hyperplane Furstenberg sets.
Here $\dimH$ denotes Hausdorff dimension.

\begin{theorem}\label{thm:hausdorff}
Let $n\ge2$, $0<s\le n-1$ and $0<t\le n$. Suppose $E\subset\R^n$
is Borel and $\HH$ is a set of affine $(n-1)$-dimensional hyperplanes satisfying
\[
\dimH\HH\ge t,\qquad \dimH(E\cap H)\ge s\quad(H\in\HH).
\]
Then
\begin{equation}\label{eq:main}
\dimH E\ge\F(s,t;n,n-1).
\end{equation}
Here $\F(s,t;n,n-1)$ is defined in \eqref{eq:index} below.
\end{theorem}

To define the index, let $\lceil s\rceil$ denote the smallest integer
greater than or equal to $s$, and set
\[
j=n-1-\lceil s\rceil\in\{0,1,\ldots,n-2\}.
\]
Define
\begin{equation}\label{eq:index}
\F(s,t;n,n-1)=
\begin{cases}
s,&0\le t\le j,\\[3pt]
\displaystyle\min\left\{s+t-j,\frac{3s+t-n+2}{2},s+1\right\},
&j<t\le n.
\end{cases}
\end{equation}

This is the hyperplane case of Gan's numerical Furstenberg
index~\cite[Definition~2.1]{Gan}, introduced for the Furstenberg
problem over prime fields. Here $n$ denotes the dimension of the
Euclidean space, while $n-1$ denotes the dimension of the hyperplanes.

When $t>j$, the integer $j$ is the number of radial reductions
needed to reach the range of the Cs\"ornyei--Stull theorem
\cite[Theorem~1.3]{CS}.
In particular, $j=0$ when $n-2<s\le n-1$, so \eqref{eq:index}
is exactly their expression in that range. For $n=2$ it is the
Ren--Wang bound~\cite[Theorem~1.1]{RW}. When $t\le j$, a $t$-dimensional family of hyperplanes can
contain the same set $E$ of dimension $s$, showing that
the bound $s$ is sharp.

The main new ingredient is the radial projection reduction in
Theorem~\ref{thm:reduction}, for $0<s\le n-2$ and $t>1$.
Each step lowers the dimension and the hyperplane-family exponent $t$
by one while preserving $s$. Combining this reduction with the
Cs\"ornyei--Stull incidence theorem
\cite[Theorem~1.1]{CS} gives the bound above.

\introsubsection{Background.}

The higher-dimensional problem builds on the planar Furstenberg set
problem, which goes back to work of Furstenberg and was presented
in Wolff's 1999 survey~\cite{Wolff}.
In its classical form, one fixes $0<s\le1$ and asks how small
the Hausdorff dimension of a set $E\subset\R^2$ can be if, in
every direction, there is a line $\ell$ with
$\dimH(E\cap\ell)\ge s$. More generally, an \((s,t)\)-Furstenberg set is defined as follows:

\begin{definition}[$(s,t)$-Furstenberg set]
For $0<s\le1$ and $0<t\le2$,
a set $E\subset\R^2$ is an $(s,t)$-Furstenberg set if there is a
family $\mathcal L$ of affine lines such that
\[
\dimH\mathcal L\ge t,\qquad
\dimH(E\cap\ell)\ge s\quad(\ell\in\mathcal L).
\]
\end{definition}

Following Wolff's work, many partial results brought the planar
Furstenberg problem closer to resolution. Orponen and
Shmerkin~\cite[Theorem~1.18]{OS} proved the sharp bound when the
associated line family has Hausdorff and packing dimensions both
equal to $t$.

They use the $ABC$ sum-product theorem to obtain a projection result,
which ensures that almost AD-regular sets have large projections
in many directions.
A multiscale argument helps to find suitable blocks of scales,
where the projection estimates give bounds
for the number of distinct parallel tubes in each of the
directions under consideration~\cite[Section~5]{OS}.

Ren and Wang~\cite{RW} removed this regularity assumption
and proved the following sharp bound. 

\begin{theorem}[Ren--Wang, {\cite[Theorem~1.1]{RW}}]
Let $0<s\le1$ and $0<t\le2$. Every $(s,t)$-Furstenberg set
$E\subset\R^2$ satisfies
\[
\dimH E\ge\min\left\{s+t,\frac{3s+t}{2},s+1\right\}.
\]
\end{theorem}
Their proof combines the almost AD-regular case with a semi-well-spaced configuration. The new estimate of the incidences is obtained through the high--low method in the paper~\cite{GSW}. To bring these estimates together, branching functions were used, which record how covering numbers grow as the scale decreases. A decomposition of these functions identifies intervals of scales on which one of the two estimates applies~\cite[Section~6]{RW}. 
\newline

It is interesting to see that in both papers, the authors got the
results by observing phenomena when lots of tubes pass through a
point, rather than looking at tubes that each contain many points.
This is because in $\R^2$, points are dual to lines. So the
Furstenberg question can either ask how large the set of points
must be, or how large the set of lines must be. In $\R^n$, when
considering affine $(n-1)$-dimensional hyperplanes,
point--hyperplane duality also interchanges points and affine
hyperplanes while preserving incidence. It therefore turns the
hyperplane Furstenberg problem into an equivalent problem about
families of hyperplanes through points. This symmetry makes
hyperplanes a natural setting in which to extend the planar
problem to higher dimensions.

\introsubsection{Recent progress in higher dimensions.} For the higher-dimensional problem,
H\'era~\cite{Hera} obtained general lower bounds for Furstenberg-type
sets associated with affine subspaces in $\R^n$. Her proof
generalizes Wolff's classical double-counting argument in two
dimensions.

Over prime fields, Gan~\cite[Theorem~1.8]{Gan} introduced a numerical
Furstenberg index and proved that the conjectured bounds in all
higher dimensions follow from the point--line Furstenberg conjecture
in $\mathbb F_p^2$. The required planar bound has since been established
by Lewko~\cite[Corollary~12]{Lewko} and follows also from the independent
incidence theorem of Miao and Xie~\cite[Theorem~1.1]{MX}.
Together with Gan's reduction, these results give the sharp exponents
in all dimensions over prime fields.

Meanwhile, in Euclidean space, Bright and
Dhar~\cite{BD} obtained new bounds for spread Furstenberg sets,
where they impose a dimension restriction on the set of directions
of the hyperplanes.

More recently, Cs\"ornyei and Stull~\cite[Theorem~1.3]{CS} extended the planar estimate
to hyperplanes in $\R^n$ when $n-2 < s \leq n-1$.

\begin{theorem}[Cs\"ornyei--Stull, {\cite[Theorem~1.3]{CS}}]
Let $n\ge2$, $n-2<s\le n-1$, and $0<t\le n$.
Suppose $E\subset\R^n$ admits a family $\HH$ of affine hyperplanes
with $\dimH\HH\ge t$ and $\dimH(E\cap H)\ge s$ for every
$H\in\HH$. Then
\[
\dimH E\ge
\min\left\{s+t,\frac{3s+t-n+2}{2},s+1\right\}.
\]
\end{theorem}
They proved new discretized versions of Marstrand's projection
theorem and Mattila's slicing theorem, which help to reduce
the incidence problem to lower dimensions.
Their discretized incidence theorem~\cite[Theorem 1.1]{CS}
is the higher-dimensional input used in this paper.

Cohen, Shmerkin, Wang, and Zakharov~\cite{CSWZ}
proved optimal dual Furstenberg estimates for affine $k$-planes
under local spread-outness assumptions on the direction measures.
In the hyperplane case, these assumptions limit the proportion of hyperplanes through a point that approximately share a common line, within each small range of directions.
After point--hyperplane duality,
their result gives
\[
\dimH E\ge s+\min\{t,1\}.
\]
The additional restrictions allow a bound that is at least
$\F(s,t;n,n-1)$, and strictly larger in some parameter ranges.
When $s+t>n$, their bound agrees with ours, since both equal $s+1$.

\introsubsection{Proof outline.}\label{sec:outline}

We prove Theorem~\ref{thm:hausdorff} by first obtaining a bound for
discretized Furstenberg configurations at scale $\delta$ and then
passing to Hausdorff dimension.
\begin{enumerate}[label=\textbf{\arabic*.}]
\item \textbf{Radial projection reduction: Theorem~\ref{thm:reduction}.}
In Section $3.1$, we first show that, by choosing centers for radial projection away from \(E\), we can obtain a projected set supporting an \(s\)-Frostman measure for each hyperplane. The centers lie on the corresponding hyperplane, at distance at least one from \(E\). Lemma 3.4 gives the required projected energy. 

The main difficulty is to find one center that works for many hyperplanes at once. 
In Section~\ref{sec:common-center}, we find a common center \(y_0\) near a family of hyperplanes whose parameters support a \((t-1)\)-Frostman measure by double counting. Projection from this same center gives a set supporting an \(s\)-Frostman measure for every hyperplane in this $(t-1)$-dimensional family.

In Section~\ref{sec:reduced-configuration}, point–hyperplane duality turns this into a configuration inside the \((n-1)\)-dimensional hyperplane \(H_{y_0}\). We obtain a \((t-1)\)-Frostman family of points, with an \(s\)-Frostman family of \((n-2)\)-dimensional hyperplanes near each point. The covering number of the new hyperplane family is controlled by that of \(E\), so a lower bound in the reduced problem also gives a lower bound for the original set.
\item \textbf{The terminal range: Proposition~\ref{prop:terminal}.}
For $n-2<s\le n-1$, we rewrite the Cs\"ornyei--Stull incidence theorem as a weighted version by dyadic pigeonholing. We then apply
\cite[Theorem~1.1]{CS} after point--hyperplane duality to fit our setting.
The resulting covering exponent is still
\begin{equation}\label{eq:top}
\min\{s+t,(3s+t-n+2)/2,s+1\},
\end{equation}
up to an arbitrarily small loss $\eps$.

\item \textbf{Iteration: Theorem~\ref{thm:finite}.}
Let $j=n-1-\lceil s\rceil$. If $t\le j$, we can place all hyperplanes around a $\lceil s \rceil$-dimensional affine subspace. If $t>j$, repeat the
reduction $j$ times, dualizing back between steps. This gives a
configuration in $\R^{n-j}$ with exponents $s$ and $t-j$, where
\[
n-j-2<s\le n-j-1.
\]
The terminal estimate therefore applies, and its exponent equals
$\F(s,t;n,n-1)$. For every $\eps>0$, choose $\eta>0$ sufficiently
small. For sufficiently small $\delta$ with $K_s,K_t\le\delta^{-\eta}$,
the logarithmic losses are absorbed and we obtain
\[
\N_\delta(E)\ge c_\eps\delta^{-\F(s,t;n,n-1)+\eps}.
\]

\item \textbf{Hausdorff dimension: Theorem~\ref{thm:hausdorff}.}
Pigeonholing a Hausdorff cover using Hausdorff content gives one scale at which both a subfamily of hyperplanes and the covered part of \(E\) on each hyperplane have sufficiently large mass. At this scale \(\delta\), the discretized theorem gives a contradiction if \(\dim_H E\) is below the discretized bound. Taking limits gives the main bound.
\end{enumerate}

Section~\ref{sec:sharpness} constructs examples attaining the bound
for all parameters in every dimension.
Subsection~\ref{sec:three-dimensional} gives the explicit formula
in $\R^3$. It also briefly discusses the author's alternative
approach to the sharp bound in three dimensions, using the high--low method
and projection estimates for almost AD-regular point sets.

\paragraph{AI assistance.}
This work began with the three-dimensional hyperplane Furstenberg problem. The author initially tried to prove the result using orthogonal projections, following the ideas in the Orponen--Shmerkin paper on the AD-regular case in two dimensions. The overall proof structure is presented in Section 3. The author intended to find a hyperplane that intersects many other hyperplanes, while the intersections can still be treated as the lower dimensional case in the same scale. During discussions of this approach, OpenAI’s GPT-5.6 suggested replacing orthogonal projections with radial projections. Further discussions with GPT generated an argument for what became Lemma 3.4. The center for radial projection is chosen outside the unit cube containing the Furstenberg set. The author subsequently used energy estimates and double counting to establish the main reduction in Section 3 in three dimensions.

In October, after knowing of the partial results on higher-dimensional Furstenberg sets used in this paper, the author generalized the three-dimensional argument to \(\mathbb R^n\). GPT-6 assisted with this extension, mainly through computations involving the radial projection map \(\Pi_y\) and tracking the various constants throughout the paper. GPT-6 was also used to write up the sharp examples in \(\mathbb R^n\), based on H\'era's example and the author's examples in three dimensions. GPT assisted with drafting and revising the manuscript and improving the language. GPT-6 also assisted with checking the proofs. The author has checked the entire proof and takes responsibility for the mathematical arguments and the final text.

\paragraph{Acknowledgements.}
I thank my advisor Ruixiang Zhang for guiding me through the literature on Furstenberg sets and for offering many insightful comments on this project. I am especially grateful for his support throughout this project and for our valuable discussions about the use of AI and its growing impact on mathematical research.
I also thank Shengwen Gan for pointing out the recent planar prime-field
results and their consequences for the higher-dimensional bounds.

\section{Definitions and notation}\label{sec:definitions}

\subsection{Scales, neighborhoods, and covering numbers}
In $\R^n$, we write $x=(x',x_n)=(x_1,x'',x_n)$, where
$x'=(x_1,\ldots,x_{n-1})\in\R^{n-1}$ and
$x''=(x_2,\ldots,x_{n-1})\in\R^{n-2}$.
We use asymptotic notation $ \lesssim,\gtrsim,\asymp$. The notation
$A\lesssim B$ means $A\le CB$ for a universal constant $C$. $A\asymp B$ means both inequalities,
and $O(\delta)$ means a quantity of absolute value at most $C\delta$.
We use the Euclidean norm $|x|$ and open balls $B(x,r) = \{p\in \R^n: \dist(p,x)<r\}$. We also write
\[
\dist(x,A)=\inf_{a\in A}|x-a|,
\qquad N_r(A)=\{x:\dist(x,A)\le r\}.
\]
The scale is a dyadic number $\delta=2^{-k}$,
where $k\ge2$ is an integer. 
A set is $\delta$-separated if any two distinct points have distance
at least $\delta$. A dyadic $\delta$-cube in $\R^n$ is
$\delta(m+[0,1)^n)$, where $m\in\mathbb Z^n$.
For a bounded set $E$, $\N_\delta(E)$ counts the dyadic
$\delta$-cubes meeting $E$; up to dimensional constants, it also
equals the minimum covering number of radius-$\delta$ balls covering $E$.
For a hyperplane family $\mathcal T$, write
$\N_\delta(\mathcal T)$ for its covering number in
hyperplane parameter space, with the metric defined below.
Thus $N_r(A)$ is a neighborhood, while $\N_\delta(E)$ is a number.
For finite collections, $\#A$ or $|A|$ denotes cardinality; for measurable subsets of \(\mathbb R^n\), \(|A|\) denotes Lebesgue measure.

Let $\mathcal D_r(A)$ denote the
half-open dyadic cubes of side length $r$ meeting $A$, and write
$\mathcal D_r=\mathcal D_r(\R^n)$. A nonempty finite collection
$\mathcal Q$ of $\delta$-cubes is a $(\delta,\alpha,K)$-set if
\[
\#\{Q\in\mathcal Q:Q\subset U\}
\le K r^\alpha\#\mathcal Q
\qquad(U \in \mathcal D_r(\R^n), \ \delta\le r\le1, r \text{ dyadic}).
\]

\subsection{Hyperplanes, their parameters, and duality}
For hyperplanes with bounded slopes, use the representation
\begin{equation}\label{eq:graph}
H_z=\{x \in \R^n:x_n=a\cdot x'+b\},\qquad
z=(a,b)\in\R^{n-1}\times\R,
\end{equation}
where $a=(a_1,a'')$ and $a,b$ stay bounded. Each hyperplane has
dimension $n-1$, while its parameter $z$ has $n$ coordinates.
In particular,
\[
z=(a,b),\qquad a=(a_1,a''),\qquad
a''=(a_2,\ldots,a_{n-1})\in\R^{n-2}.
\]
We may measure distance between hyperplanes in these coordinates by
\[
\operatorname{dist}_{\mathrm{hyp}}(H_z,H_{\widetilde z})
=|z-\widetilde z|
=\sqrt{|a-\widetilde a|^2+|b-\widetilde b|^2}.
\]
In particular, a set of hyperplanes is $\delta$-separated if the distance between any two distinct hyperplanes is at least $\delta$.
In the paper, there will be three related objects: the exact hyperplane $H_z$; its
$\delta$-thickening $N_{\delta}(H_z)$; and a dyadic
$\delta$-hyperplane, meaning the family of hyperplanes whose
parameters lie in one dyadic $\delta$-cube, as in \cite{CS}.
The last can be seen as a cell in parameter space. 

A set of dyadic $\delta$-hyperplanes in $\R^n$
is called a $(\delta,s,K)$-set if the corresponding set in parameter space is a
$(\delta,s,K)$-set.
The dyadic cube condition and the ball condition in \(\eqref{eq:frostman}\) are actually equivalent after multiplying the Frostman constant \(K\) by a constant depending only on the dimension. A dyadic \(\delta\)-hyperplane intersects a cube if at least one exact hyperplane in the family intersects that cube.

Point--hyperplane duality sends $H_z$ to the point $z$, and sends
a point $x$ to the dual hyperplane
\[
x^*=\{(a,b):a\cdot x'+b=x_n\}.
\]
For a viewpoint \(y_0=(u,v)\in \R^{n-1} \times \R\), we use \(H_{y_0}\) to denote its dual hyperplane in parameter space:
\[
H_{y_0}:=y_0^*=\{(a,b):a\cdot u+b=v\}.
\]
We write
\[
y_0=(u,v)=(u_1,u'',v),\qquad
u''=(u_2,\ldots,u_{n-1})\in\R^{n-2}.
\]
Its basic property is $x\in H_z$ if and only if $z\in x^*$.

Additionally, points of $x^*\cap H_{y_0}$ parameterize the original hyperplanes
passing through both $x$ and $y_0$.

\subsection{Frostman measures and weighted $(s,t)$-configurations}
All measures called probabilities have total mass $1$, and $\supp\mu$ denotes the support of $\mu$. For $\mu(A)>0$, restricting to $A$ and dividing by
its mass gives another probability measure on $A$,
$\mu_A(B)=\mu(A\cap B)/\mu(A)$.
For a measurable map $f$, the pushforward is
\[
(f_\#\mu)(B)=\mu(f^{-1}(B)).
\]

\begin{definition}[Frostman measure at scale $\delta$]
Fix $0<\delta<1/2$, an exponent $\alpha\ge0$, and a constant $K\ge1$.
A Borel probability measure $\mu$ on $\R^n$ is a Frostman
$(\delta,\alpha,K)$ measure if, for every center $x\in\R^n$,
\begin{equation}\label{eq:frostman}
\mu(B(x,r))\le K r^\alpha\qquad(\delta\le r\le1).
\end{equation}
Throughout this paper, a Frostman measure has total mass
$\mu(\R^n)=1$.

For a finite family of balls $B=\{B(z,\delta)\}_z$, with repetitions allowed,
choose weights $m_z\ge0$ with $\sum_z m_z=1$, and distribute each weight
uniformly over its ball. This defines the probability measure
\[
\mu(A)=\sum_z m_z\,
\frac{|A\cap B(z,\delta)|}{|B(z,\delta)|}
\qquad(A\subset\R^n).
\]
Here $|\cdot|$ is Lebesgue measure, and the contributions add where
balls overlap. This measure is a Frostman $(\delta,\alpha,K)$
measure when it satisfies \eqref{eq:frostman}.

\end{definition}

\begin{definition}[Weighted $(s,t)$ configuration]
\label{def:configuration}
A weighted $(s,t)$ configuration $(E,H)$ at scale $\delta$ consists
of a bounded $\delta$-separated point set $E\subset\R^n$ and a $\delta$-separated hyperplane family
$H \subset \mathcal A(n,n-1)$ in $\mathbb{R}^n$ with bounded parameters. 
They satisfy the following two conditions: 
\begin{itemize}
\item The parameter set of hyperplanes $Z = \{z: H_z \in H\}$ supports a Frostman $(\delta, t, K_t)$ measure $\nu$.

\item For every $H_z\in H$, the subset $E\cap N_\delta(H_z)$ supports a Frostman $(\delta, s, K_s)$ measure $\mu_z$.
\end{itemize}
\end{definition}

Neither measure is required to be a uniform counting measure;
this flexibility is necessary after projection.
We use the word “weighted” because different points and hyperplane parameters may have different masses.

Later, we always use $\nu$ to describe the $t$-dimensional measure of the hyperplanes, while
each $\mu_z$ describes the $s$-dimensional probability measure of points near its
hyperplane $H_z$.

\begin{definition}[Dual weighted $(s,t)$ configuration]
\label{def:dualconfiguration}
In $\R^m$, a dual weighted $(s,t)$ configuration $(P,\mathcal T)$ at scale $\delta$
consists of a bounded $\delta$-separated point set $P$ and a $\delta$-separated
hyperplane family $\mathcal T$ of bounded parameters. The points support a Frostman
$(\delta,t,K_t)$ measure. At each $p\in P$, there is a subset of $\mathcal T$ intersecting $B(p,\delta)$, which supports a Frostman
$(\delta,s,K_s)$ measure in hyperplane parameter space.
\end{definition}

\textit{Remark:} Notice that the two forms of weighted configuration are the same after duality.
Meanwhile, another formulation of a weighted $(s,t)$-configuration replaces $x \in E$ by $B(x,\delta)$
and replaces $z \in Z$ by the parameter ball $B(z,\delta)$.
A dyadic $\delta$-hyperplane is instead represented by a dyadic parameter cube.
The first condition is unchanged. Define
$
\mathcal B=\{B(x,\delta):x\in E\}.
$
The second condition requires, for each $H_z$, that the family
\[
\mathcal B_z=\{B\in\mathcal B:B\cap N_{\delta}(H_z)\ne\varnothing\}
\]
form a Frostman $(\delta,s,K_s)$ measure, with its own weights
counted as above. The two different formulations are basically equivalent since moving mass between each ball
and its center changes the constant $K_\alpha$ in the Frostman measure by at most a constant factor. 

\begin{definition}[Truncated energy]
For a probability measure $\mu$, an exponent $\alpha\ge0$, and
$\delta>0$, define the truncated energy by
\[
I_\alpha^\delta(\mu)
=\iint (|x-\tilde x|+\delta)^{-\alpha}
\dd\mu(x)\dd\mu(\tilde x).
\]
If $\mu$ is finitely supported, with weights $m_x=\mu(\{x\})$,
then
\[
I_\alpha^\delta(\mu)
=\sum_{x,\tilde{x}\in\supp\mu}
m_xm_{\tilde x}(|x-\tilde{x}|+\delta)^{-\alpha}.
\]
\end{definition}

\subsection{Radial projections}
The unit-direction radial projection from a viewpoint $y$ is
\[
\pi_y(x)=\frac{x-y}{|x-y|}\in S^{n-1},\qquad x\ne y.
\]
For $x\in[-1,1]^n$ and a viewpoint $y = (y_1,y') \in \R^n$ with $3\le y_1\le4$,
we can also represent this projection in the affine hyperplane
$L=\{x_1=0\}\cong\R^{n-1}$. Since $x_1\ne y_1$, the line
through $y$ and $x$ meets $L$ at a unique point $(0,\Pi_y(x)) \in L$. We denote its
last $n-1$ coordinates by $\Pi_y(x)$; these coordinates represent
the same direction as $\pi_y(x)$. 

The trace of one hyperplane on another means their intersection. Fix a viewpoint $y_0$.
For each \(x\in E\), the line through \(y_0\) and \(x\) meets \(L=\{x_1=0\}\) at a unique point \((0,w_x)\), where \(w_x=\Pi_{y_0}(x)\) records its last \(n-1\) coordinates,
\[
(0,w_x)=(0,w_x'',w_{x,n})\in L,\qquad
w_x=\Pi_{y_0}(x),\qquad
w_x''=(w_{x,2},\ldots,w_{x,n-1}).
\]

The double-prime vectors always have $n-2$ coordinates in this notation. Note that the trace $x^*\cap H_{y_0}$ is an $(n-2)$-dimensional hyperplane
inside $H_{y_0}$, and can be represented by
\[ x^*\cap H_{y_0}
=\{(a'',b) \in \R^{n-2} \times \R:b=w_{x,n}-w_x''\cdot a''\}.
\]
In the trace equation $b=w_{x,n}-w_x''\cdot a''$, both vectors
in the dot product belong to $\R^{n-2}$.

This is the main reason we chose $\Pi_y(x)$ to represent the directions instead of using $\pi_y(x)$ directly. This parameterization is also shown in the
reduction argument.

Meanwhile, the use of radial projection in the main proof is inspired by the fact that the dimension of the radial projection of an $s$-dimensional point set is almost always equal to the dimension of the original set. More specifically,
\begin{theorem}[Radial projections in $\R^n$ {\cite[Theorem 2]{BG}}]
\label{thm:radialdimension}
Let $A\subset\R^n$ be Borel, with $0<s=\dimH A\le n-1$. Then
\[
\dimH\{y\in\R^n\setminus A:\dimH\pi_y(A)<s\}
\le\lceil s\rceil.
\]
\end{theorem}
In particular, if \(A\subset H_z\) has Hausdorff dimension \(0<s\le n-2\), then radial projection preserves its dimension for almost every viewpoint on \(H_z\) outside \(A\):
$
\dim_{\mathrm H}\pi_y(A)=s.
$
The exceptional viewpoints form a set of dimension at most \(n-2\), so they have zero \((n-1)\)-dimensional measure on \(H_z\).

Throughout the proof in the next section, we choose viewpoints a fixed positive distance from the original bounded set. We also use the word “center” to mean a projection viewpoint.
\section{A discretized reduction theorem}

The main ingredient of Theorem~\ref{thm:hausdorff} is the following result. After converting the Furstenberg set into a discretized $\delta$-set, we use this result to extract a lower-dimensional dual Furstenberg set. The lower-dimensional Furstenberg bound then gives a bound for the dimension of the original set. 
\begin{theorem}\label{thm:reduction}
Let $n\ge3$, $0<s\le n-2$ and $1<t\le n$. Given a weighted $(s,t)$
configuration $(E,H)$ at scale $\delta$ in $\R^n$, there exist a viewpoint $y_0$ and a dual weighted $(s,t-1)$
configuration $(P,\mathcal T)$ at scale $\kappa\delta$ inside the $(n-1)$-dimensional hyperplane
$H_{y_0}\cong\R^{n-1}$. Moreover, the family $(P, \mathcal T)$ satisfies
\begin{equation}\label{eq:losses}
\begin{aligned}
K_s'&\lesssim K_s\log^{2}(1/\delta),\qquad K_t'\lesssim K_t\log(1/\delta), \qquad
\N_{\kappa\delta}(\mathcal T)&\lesssim\N_\delta(E).
\end{aligned}
\end{equation}
Note that the implied constants and \(\kappa\) depend only on \(s,t,n\) and the fixed bounded regions containing \(E\) and the hyperplane
parameter set \(Z\). They are independent of \(\delta\).
The last inequality also means that
\[
\#\mathcal T\lesssim\#E.
\]

\end{theorem}
We prove this theorem in the remainder of this section.

\noindent\textit{Proof of Theorem~\ref{thm:reduction}.}
Normalize $E\subset[-1,1]^n$ and use the notation in \eqref{eq:graph}.
By Definition~\ref{def:configuration}, both $E$ and the parameter
set $Z$ are already $\delta$-separated.
Denote $x = (x_1,x'',x_n), a = (a_1, a'')$, while $x' = (x_1, x'')$. The parameter space for the radial projection is $L = \{x_1=0\}\cong\R^{n-1}$.

Set $$U = [3,4]\times [-1,1]^{n-2},$$ which will be the region where we take viewpoints' first $n-1$ coordinates from. 

For each hyperplane $H_z = \{x: x_n = a \cdot x' + b\}$ in $\mathbb R^n$, 
take a viewpoint $y$ in $H_z$ such that $y = (u, a \cdot u + b)$ for some $u \in U$.

Given this viewpoint $y$, for a point $x \in E \cap N_\delta(H_z)$, 
the line from $y$ to $x$ meets $L$ at a unique point 
since the first coordinate of $y$ is away from the unit cube. 
As stated in Section $2$, we use $(0,\Pi_y(x))$ to denote the point on $L$ 
that is on the line $\{y+t(x-y): t \in \mathbb R\}$, so
\begin{equation}\label{eq:projection}
\Pi_y(x)=\left(
\frac{u_1x''-x_1u''}{u_1-x_1},\,
\frac{u_1x_n-x_1v}{u_1-x_1}
\right).
\end{equation}
Clearly, if \(x\in H_z\), then \((0,\Pi_y(x))\in H_z\cap L\).
Because $|u_1-x_1|\ge2$ for any pair of $y$ and $x \in E\cap H_z$, the maps are uniformly Lipschitz in both $x$ and viewpoints $y$, and have bounded images. The values change by $O(\delta)$ when
the point $x$ changes by $O(\delta)$ within the neighborhood $N_{\delta}(H_z)$.

Thus, for $y\in H_z$, the original pushforward $(\Pi_y)_\#\mu_z$ has support on 
an $O(\delta)$-neighborhood of the hyperplane $H_z\cap L$ in $L$.

We next show that distances between flat projections $\Pi_y(x)$ are comparable to distances between the corresponding unit directions $\pi_y(x)$.
The following lemma also relates the projected distance to the line joining the original points.
\begin{lemma}\label{lem:geometry}
Fix a constant $C>0$. For distinct $x,\widetilde x\in[-1,1]^n$, let $\ell$ be their joining
line. Uniformly for the bounded viewpoints $y=(y_1,y'',y_n) \in [3,4] \times [-1,1]^{n-2}\times[-C,C]$,
\begin{equation}\label{eq:distancecompare}
|\Pi_y(x)-\Pi_y(\widetilde x)|
\asymp |x-\widetilde x|\dist(y,\ell).
\end{equation}
The points $x,\widetilde x$ need not lie on $H_z$. In particular, $|\Pi_y(x)-\Pi_y(\widetilde x)|\le C'|x-\widetilde x|$. 
The comparison constants and $C'$ depend only on $n$ and $C$.
\end{lemma}
\begin{proof}
Set the unit directions
$
\pi_y(x) = \omega=\frac{x-y}{|x-y|},
\pi_y(\widetilde x) = \widetilde\omega=\frac{\widetilde x-y}{|\widetilde x-y|}.
$
The bounded domains give $2\le|x-y|,|\widetilde x-y|\le C+2n+2$
and $|x_1-y_1|,|\widetilde x_1-y_1|\ge2$. Hence
$|\omega_1|,|\widetilde\omega_1|\ge \frac{2}{C+2n+2}$;
these two first coordinates have the same sign.
Notice that by the formula for \(\Pi_y\),
\[
|\Pi_y(x)-\Pi_y(\widetilde x)|
=
y_1\left|
\frac{\omega}{\omega_1}
-\frac{\widetilde\omega}{\widetilde\omega_1}
\right|.
\]
These bounds give the following upper and lower estimates:
\begin{equation}\label{eq:flatdirection}
\frac{y_1}{2}|\omega-\widetilde\omega|
\le |\Pi_y(x)-\Pi_y(\widetilde x)|
\le \frac{y_1(C+2n+2)^2}{4}|\omega-\widetilde\omega|.
\end{equation}
Thus the flat and unit-direction distances are uniformly comparable as $y_1$ is always inside $[3,4]$.
If $\theta$ is the angle between $\omega$ and $\widetilde\omega$, then
\begin{equation}\label{eq:anglecompare}
\sin\theta\asymp|\omega-\widetilde\omega|.
\end{equation}
Computing the area of the triangle $y,x,\widetilde x$ in two ways gives
\begin{equation}\label{eq:triangle}
\sin\theta
=\frac{|x-\widetilde x|\dist(y,\ell)}
{|x-y|\,|\widetilde x-y|}.
\end{equation}
Since $3\le y_1\le4$ and both denominator factors lie in $[2,C+2n+2]$,
we conclude
\[
|\Pi_y(x)-\Pi_y(\widetilde x)|
\asymp|\omega-\widetilde\omega|
\asymp\sin\theta
\asymp|x-\widetilde x|\dist(y,\ell).
\]
All constants depend only on the fixed bounded domains $E$ and $U$.
\end{proof}

We will use this distance comparison to estimate the energy of the projected measures. Before doing so, we recall the connection between the Frostman condition and truncated energy, given by the following lemma.
Recall that the truncated energy is defined to be \[
I_\alpha^\delta(\mu)=\iint(|x-x'|+\delta)^{-\alpha}
\dd\mu(x)\dd\mu(x').
\]

\begin{lemma}\label{lem:energy}
On a fixed bounded region:
\begin{enumerate}[label=(\roman*)]
\item A Frostman $(\delta,\alpha,K)$ measure satisfies
$I_\alpha^\delta(\mu)\le C K \log(1/\delta)$.
\item If $I_\alpha^\delta(\mu)\le A$, we can restrict $\mu$ to a subset
of mass at least $1/2$ and normalize to obtain a
Frostman $(\delta,\alpha,CA)$ measure.
\end{enumerate}
\end{lemma}
\begin{proof}
For (i), sum the kernel over annuli of radii
$\delta,2\delta,4\delta,\ldots,1$, using \eqref{eq:frostman}; each
annulus contributes at most $CK$. Distances greater than $1$ contribute
at most a constant. There are $O(\log(1/\delta))$ annuli.

For (ii), put $C(x)=\int(|x-x'|+\delta)^{-\alpha}\dd\mu(x')$.
Markov's inequality gives $\mu(G)\ge1/2$ for $G=\{x: C(x)\le2A\}$. Restrict $\mu$ to $G$ and normalize its mass.
For a ball $B(c,r)$ meeting $G$, choose $x\in G\cap B(c,r)$.
Then $(3r)^{-\alpha}\mu(G\cap B(c,r)) \le C(x) \le 2A$ for $r\ge\delta$.
\end{proof}

\proofsubsection{Good projection centers for each hyperplane.}\label{sec:good-projection-centers}

The following lemma is the first key lemma in the proof of Theorem~\ref{thm:reduction}.
It bounds the average truncated $s$-energy of the projected measures
$(\Pi_y)_\#\mu_z$ by $CK_s\log^2(1/\delta)$.
Together with Lemma~\ref{lem:energy}(ii), this gives a Frostman
$(\delta,s,CK_s\log^2(1/\delta))$ measure on a subset of the projected
set for some viewpoint.

Partition $U= [3,4] \times [-1,1]^{n-2}$ into dyadic cubes of side $\delta$, and let $\mathcal U_\delta$ be their centers.
This is a $\delta$-separated grid with
$
\#\mathcal U_\delta=2^{n-2}\delta^{-(n-1)}.
$
Given a hyperplane $H_z$, each $u\in\mathcal U_\delta$
determines a viewpoint $y_z(u)=(u,a\cdot u+b)\in H_z$. So there are $\#\mathcal U_\delta$ many points on
$H_z$ that can be treated as viewpoints.
\begin{lemma}\label{lem:radial}
Let $n\ge3$ and $0<s\le n-2$, and let \(Z\subset\mathbb R^n\) be a fixed bounded set of hyperplane parameters. For each \(z\in Z\), suppose \(\mu_z\) is a Frostman \((\delta,s,K_s)\) measure supported on a \(\delta\)-separated set inside $N_\delta(H_z)$,
\[
E_z\subset[-1,1]^n\cap N_{\delta}(H_z).
\]
With $U$, $\mathcal U_\delta$, and the viewpoints $y_z(u)$
defined above, we have
\begin{equation}\label{eq:radialaverage}
\frac{1}{\#\mathcal U_\delta}
\sum_{u\in\mathcal U_\delta}
I_s^\delta\bigl((\Pi_{y_z(u)})_\#\mu_z\bigr)
\le CK_s\log^2(1/\delta).
\end{equation}
If $s<n-2$, the right-hand side can be replaced by
$CK_s\log(1/\delta)$. The constant $C$ depends only on
$n,s$ and the bound for hyperplane parameters.
\end{lemma}
\begin{proof}
This proof uses double counting. Fix distinct points $x,\widetilde x\in E_z$, so
$|x-\widetilde x|\ge\delta$. For $r\ge\delta$,
Lemma~\ref{lem:geometry} gives
\begin{equation}\label{eq:tube}
\begin{aligned}
&\#\{u\in\mathcal U_\delta:
|\Pi_{y_z(u)}(x)-\Pi_{y_z(u)}(\widetilde x)|\le r\}
\\
&\qquad\lesssim \delta^{-(n-1)}
\left(\frac{r}{|x-\widetilde x|}\right)^{n-2},
\end{aligned}
\end{equation}
where we count the viewpoints \(y_z(u)\), \(u\in\mathcal U_\delta\), lying in a tube of radius \(\frac{r}{|x-\widetilde x|}\) around the line joining \(x\) and \(\widetilde x\).
When $r\ge|x-\widetilde x|$, this bound follows simply from
$\#\mathcal U_\delta\asymp\delta^{-(n-1)}$.

Group the truncated projected distances
$|\Pi_{y_z(u)}(x)-\Pi_{y_z(u)}(\widetilde x)|+\delta$
into dyadic intervals $[r,2r)$. The above bound holds for every distance $r$
between $\delta$ and $C'|x-\widetilde x|$ (where $C'$ is the constant from Lemma~\ref{lem:geometry}). Counting these viewpoints with \eqref{eq:tube} gives
\[
\begin{aligned}
&\frac{1}{\#\mathcal U_\delta}
\sum_{u\in\mathcal U_\delta}
\bigl(|\Pi_{y_z(u)}(x)-\Pi_{y_z(u)}(\widetilde x)|+\delta\bigr)^{-s}
\\
&\quad\le\frac{1}{\#\mathcal U_\delta}
\sum_{\substack{r\ \mathrm{dyadic}\\\delta\le r\le C'|x-\widetilde x|}}
r^{-s}
\#\{u \in \mathcal U_\delta: 
|\Pi_{y_z(u)}(x)-\Pi_{y_z(u)}(\widetilde x)|\le2r\}
\\
&\quad\lesssim |x-\widetilde x|^{-(n-2)}
\sum_{\substack{r\ \mathrm{dyadic}\\\delta\le r\le C'|x-\widetilde x|}}
r^{n-2-s}
\\
&\quad\lesssim
\begin{cases}
|x-\widetilde x|^{-s},&s<n-2,\\
\log(1/\delta) |x-\widetilde x|^{-s},&s=n-2.
\end{cases}
\end{aligned}
\]
The averaged truncated energy is the weighted sum over all pairs $(x,\widetilde x)$:
\[
\begin{aligned}
&\frac{1}{\#\mathcal U_\delta}
\sum_{u\in\mathcal U_\delta}
I_s^\delta\bigl((\Pi_{y_z(u)})_\#\mu_z\bigr)
\\
&\quad=\frac{1}{\#\mathcal U_\delta}
\sum_{u\in\mathcal U_\delta}
\sum_{x,\widetilde x\in E_z}
\frac{\mu_z(\{x\})\mu_z(\{\widetilde x\})}
{\bigl(|\Pi_{y_z(u)}(x)-\Pi_{y_z(u)}(\widetilde x)|+\delta\bigr)^s}
\\
&\quad\lesssim \log(1/\delta)
\sum_{x,\widetilde x\in E_z}
\frac{\mu_z(\{x\})\mu_z(\{\widetilde x\})}
{\bigl(|x-\widetilde x|+\delta\bigr)^s}
\\
&\quad = \log(1/\delta)I_s^\delta(\mu_z)
\\
&\quad\lesssim K_s\log^2(1/\delta).
\end{aligned}
\]
The last inequality follows from the Frostman bound for $\mu_z$
and Lemma~\ref{lem:energy}(i). Diagonal pairs are included:
when $x=\widetilde x$, both truncated kernels equal $\delta^{-s}$.
If $s<n-2$, the preceding pairwise estimate in the end gives $CK_s\log(1/\delta)$.
\end{proof}

By the finite-sum form of Markov's inequality, there is a set in
$\mathcal U_\delta$ of at least half the points where the energy is at most $2CK_s\log^{2}(1/\delta)$. 
Denote $
\mathcal U_{\delta,z}
=
\left\{
u\in\mathcal U_\delta:
I_s^\delta\bigl((\Pi_{y_z(u)})_\#\mu_z\bigr)
\le 2CK_s\log^2(1/\delta)
\right\}.
$
Inside the $(s,t)$-configuration,
choose a fixed $c>0$ sufficiently small that $c\le1/4$ and
$c\sqrt{1+|a|^2}\le1$ for every parameter $z = (a,b)$ of hyperplanes.
Define the good-center set directly by balls around the selected viewpoints:
\begin{equation}\label{eq:goodcenters}
G_z=\bigcup_{u\in\mathcal U_{\delta,z}}
B\bigl((u,a\cdot u+b),c\delta\bigr) \subset N_{\delta}(H_z).
\end{equation}
The centers are $\delta$-separated, so these balls are disjoint and
\[
|G_z|\asymp\delta^n\#\mathcal U_{\delta,z}\asymp\delta.
\]
The choice of $c$ ensures that every $y=(y',y_n)\in G_z$ satisfies
$y'\in U$ and $|y_n-a\cdot y'-b|\le\delta$. Thus $y \in N_\delta(H_z)$.
Every point of $G_z$ is within $c\delta$ of a selected viewpoint,
so the $O(\delta)$ stability of the projection gives, for every $y\in G_z$,
\begin{equation}\label{eq:goodenergy}
I_s^\delta((\Pi_y)_\#\mu_z)\le CK_s\log^{2}(1/\delta).
\end{equation}

\proofsubsection{Selecting one good center for the whole family.}\label{sec:common-center}
After finding the good center set for each hyperplane $H_z$, we next find one good center that works for a large subfamily of hyperplanes. The argument uses the fact that two hyperplanes in $\mathbb{R}^n$ don't overlap too much. Below, $|\cdot|$ denotes Lebesgue measure in $\mathbb{R}^n$. \begin{lemma}[Two thickened planes via duality]\label{lem:overlap}
For any two good center sets inside distinct hyperplanes, the volume of the intersection satisfies
\begin{equation}\label{eq:overlap}
|G_z\cap G_{z'}|\lesssim\frac{\delta^2}{|z-z'|+\delta}.
\end{equation}
\end{lemma}
\begin{proof}
Let $\mathcal Q$ be the set of $\delta$-cubes meeting $G_z\cap G_{z'}$,
with centers $y_Q=(u_Q,v_Q)$. Under point--hyperplane duality,
two hyperplanes $H_z,H_{z'}$ that contain the good center sets become the points $z,z'$, while
\[
y_Q\longleftrightarrow
 y_Q^*=\{(a,b):b=v_Q-u_Q\cdot a\}.
\]
$y_Q^*$ meets both $B(z,C\delta)$ and $B(z',C\delta)$.
The two-ball hyperplane count, followed by multiplying by the
$\delta$-cube volume $\delta^n$, gives
\[
|G_z\cap G_{z'}|
\le\delta^n\#\mathcal Q
\lesssim\delta^n\frac{1}{(|z-z'|+\delta)\delta^{n-2}}
=\frac{\delta^2}{|z-z'|+\delta}.
\qedhere
\]
\end{proof}
Next, we find one center $y_0$ such that the hyperplanes whose good-center sets contain \(y_0\) look like a $(t-1)$-dimensional set. Again, truncated energy will be crucial for finding the viewpoint we want. 

Recall that $\nu$ is a measure on the hyperplane parameter space $Z$ and is assumed to be a Frostman $(\delta, t, K_t)$ measure. 
For each viewpoint $y$, define the total $\nu$-mass of hyperplanes
for which $y$ is a good center, together with the corresponding
unnormalized truncated $(t-1)$-energy:
\begin{align*}
M(y)&=\int \mathbf1_{G_z}(y)\dd\nu(z),\\
\mathcal E(y)&=\iint \mathbf1_{G_z}(y)\mathbf1_{G_{z'}}(y)
(|z-z'|+\delta)^{-(t-1)}\dd\nu(z)\dd\nu(z').
\end{align*}

Write $Z_y=\{z\in Z:y\in G_z\}$. For $M(y)>0$, let $\nu_y$
be the restriction of the measure $\nu$ to $Z_y$, divided by $M(y)$.

We expect the normalized energy of some $\nu_y$ to be bounded by a constant apart from a logarithmic factor. Then it implies that the family of hyperplanes lying within distance \(\delta\) of \(y\) supports a Frostman \((\delta,t-1,CK_t\log(1/\delta))\) measure. Choose a fixed constant \(D>0\), large enough that every good-center set \(G_z\) lies in the box
\begin{equation}\label{eq:omega}
\Omega:=U\times[-D,D]
=[3,4]\times[-1,1]^{n-2}\times[-D,D].
\end{equation}
All integrals over \(\Omega\) below use \(n\)-dimensional Lebesgue measure. 

Fubini and Cauchy--Schwarz give
\begin{equation}\label{eq:centerintegrals}
\int_\Omega M\gtrsim\delta,\qquad
\int_\Omega M^2\gtrsim\delta^2.
\end{equation}
For the total truncated energy, exchange the order of integration, then apply
Lemma~\ref{lem:overlap} and Lemma~\ref{lem:energy}(i):
\[
\begin{aligned}
\int_\Omega\mathcal E(y)\dd y
&=\iint \frac{|G_z\cap G_{z'}|}{(|z-z'|+\delta)^{t-1}}
\dd\nu(z)\dd\nu(z')\\
&\lesssim\delta^2\iint (|z-z'|+\delta)^{-t}
\dd\nu(z)\dd\nu(z')\\
&=\delta^2 I_t^\delta(\nu)
\lesssim K_t\log(1/\delta)\delta^2.
\end{aligned}
\]
Notice that the normalized truncated energy is $$ I_{t-1}^{\delta}(\nu_y) = \frac{1}{M(y)^2}\int\int_{Z_y \times Z_y} (|z-z'| + \delta)^{-(t-1)}d\nu(z)d\nu(z') = \frac{\mathcal E(y)}{M(y)^2}. $$
Combining the bounds for $\int_\Omega M^2$ and
$\int_\Omega\mathcal E$, we find at least one $y\in\Omega$ with $M(y)>0$ such that 
\begin{equation}\label{eq:ratio}
I_{t-1}^{\delta}(\nu_y) = \mathcal E(y)/M(y)^2\le CK_t\log(1/\delta).
\end{equation}

Fix one such viewpoint and denote it by $y_0\in\Omega$.
The probability measure
\[
\dd\nu_{y_0}(z)=M(y_0)^{-1}\mathbf1_{G_z}(y_0)\dd\nu(z)
\]
has truncated $(t-1)$-energy bounded by \eqref{eq:ratio}.
This is the common center used in the rest of the proof.

\proofsubsection{The configuration in lower dimension by duality.}\label{sec:reduced-configuration}

We now return to the proof of Theorem~\ref{thm:reduction} and complete it using the center $y_0$ found in Section~\ref{sec:common-center}.
Write $y_0=(u,v)\in\R^n$, with $u=(u_1,u'')\in\R^{n-1}$.
We construct the dual weighted $(s,t-1)$ configuration of
Definition~\ref{def:dualconfiguration}: a $(t-1)$-Frostman family
of balls  inside $H_{y_0}$ and, near each ball, an $s$-Frostman family of $(n-2)$-dimensional
hyperplanes. We break down the following argument into five steps.

\proofstep*{\textbf{Step 1:} The $(t-1)$-Frostman point family}
Under point–hyperplane duality, each original hyperplane \(H_z\) becomes its parameter point \(z\), while each 
point $x\in E$ becomes the hyperplane
\[
x^*=\{(a,b):b=x_n-a\cdot x'\},
\qquad
H_{y_0}=\{(a,b):b=v-u\cdot a\}.
\]
Put
\[
Z_1=\{z\in Z:y_0\in G_z\}.
\]
By \eqref{eq:goodcenters}, for each $z \in Z_1$, $B(z,\delta)$ meets $H_{y_0}$ and
$
B(z,\delta)\subset N_{2\delta}(H_{y_0}).
$
The choice of \(y_0\) gives
\[
I_{t-1}^{\delta}(\nu_{y_0})
\lesssim K_t\log(1/\delta).
\]
By Lemma 3.3(ii), restricting and normalizing this measure gives a Frostman
$
(\delta,t-1,CK_t\log(1/\delta))
$
measure \(\lambda\) on a subset \(Z_2\subset Z_1\).
We now project each \(z\in Z_2\) to its nearest point \(p_z\in H_{y_0}\). Since
$
|p_z-z|\lesssim\delta,
$
transferring the weights to these points preserves the Frostman bound up to a fixed factor. We keep the original label $Z_2$ when referring to points on $H_{y_0}$ and keep \(z\) when referring to \(\mu_z\).

\proofstep*{\textbf{Step 2:} How hyperplanes intersect with $H_{y_0}$}

We next describe the hyperplanes that will pass near these
points. On $H_{y_0}$, the first coordinate is determined by
the others:
\[
a_1=\frac{v-u''\cdot a''-b}{u_1}.
\]
Thus we can use $(a'',b)\in\R^{n-1}$ as coordinates on $H_{y_0}$ since $u \in U = [3,4]\times [-1,1]^{n-2}$.
These coordinates preserve distances up to fixed factors.

Next, we record the points on $x^*\cap H_{y_0}$ by using the last $n-1$ coordinates. For each $x\in E$, substituting this expression into the equation
of $x^*$, together with the equation of $H_{y_0}$, gives the coordinate representation
\[
x^*\cap H_{y_0} \rightarrow T_{w_x}:=\{(a'',b) \in \mathbb R^{n-1}:b=w_{x,n}-w_x''\cdot a''\},
\qquad w_x=\Pi_{y_0}(x),
\]
where
\[
w_x=(w_x'',w_{x,n})
=\left(\frac{u_1x''-x_1u''}{u_1-x_1},
         \frac{u_1x_n-x_1v}{u_1-x_1}\right).
\]
Thus each original point $x\in E$ produces a projected point $w_x$,
and duality produces an intersection hyperplane $T_{w_x}$ of dimension $n-2$. The projection $w_x$ records
the whole intersection hyperplane, while points $(a'',b)$ represent coordinates
of its points. Notice that different original points $x$ may produce the same $w_x$ and the same trace.

Using the coordinates \((a'',b) \in \R^{n-1}\) on \(H_{y_0}\), we measure the distance between two trace hyperplanes by the distance between their coefficients.
Again, since \(a_1\) is determined by \((a'',b)\), this loses no information and preserves distances up to fixed factors. For \(x,\widetilde x\in E\), the corresponding traces have parameters \(w_x=\Pi_{y_0}(x)\) and \(w_{\widetilde x}=\Pi_{y_0}(\widetilde x)\). Hence
\[
\operatorname{dist}_{\mathrm{hyp}}(T_{w_x},T_{w_{\widetilde x}})
\asymp |w_x-w_{\widetilde x}|
=|\Pi_{y_0}(x)-\Pi_{y_0}(\widetilde x)|.
\]

\proofstep*{\textbf{Step 3:} An $s$-Frostman family of trace hyperplanes}

Fix an original parameter $z\in Z_2$. Since $y_0\in G_z$,
the projected probability
$
\tau_z=(\Pi_{y_0})_\#\mu_z
$
satisfies
\[
I_s^\delta(\tau_z)\lesssim K_s\log^2(1/\delta).
\]
Lemma~\ref{lem:energy}(ii) gives a subset $W_z\subset\{w_x = \Pi_{y_0}(x):x\in\supp\mu_z\}
$ of the projected
set supporting a Frostman
$
(\delta,s,CK_s\log^2(1/\delta))
$
measure $\rho_z$.

Assign the weight of each $w\in W_z$ to the corresponding trace 
hyperplane $T_w$. By the distance comparison above, the trace family
$
\mathcal T_z=\{T_w:w\in W_z\}
$
has the same Frostman bound. 

Hence, the family of points
$\{p_z:z\in Z_2\}\subset H_{y_0}$ supports a Frostman
$(\delta,t-1,CK_t\log(1/\delta))$ measure.
For each $z\in Z_2$, the corresponding family $\mathcal T_z$
of trace hyperplanes supports a Frostman
$(\delta,s,CK_s\log^2(1/\delta))$ measure $\rho_z$.
It remains to check that these hyperplanes on $H_{y_0}$ meet the ball
corresponding to $p_z$.

\proofstep*{\textbf{Step 4:} The trace hyperplanes meet the $\kappa \delta$-ball}
We claim that $z$ is indeed within the $\kappa\delta$-neighborhood of $T_{w_x}$ for some constant $\kappa$
that is independent of $\delta$. 
Thus every member of $\mathcal T_z$ meets a ball of radius
comparable to $\delta$ about the point corresponding to $p_z$.

Choose \(k\ge1\), independent of \(\delta\), such that both \(\Omega\) in \eqref{eq:omega} and the original parameter set \(Z\) lie in \([-k,k]^n\). We show that every member of $\mathcal T_z$ lies in $H_{y_0}$ and meets
$B(z,3nk\delta)\cap H_{y_0}$.

This can be done by examining the distance directly. Let $z = (a,v - a \cdot u) \in H_{y_0}$ and
$x \in N_{\delta}(H_z)$. Then 
the distance between $z$ and $T_{w_x}$ is equivalent to the distance
for the last $n-1$ coordinates up to a constant, so write $\kappa=3nk$. We obtain
\[
\begin{aligned}
\operatorname{dist}\bigl(z, T_{w_x} \bigr)
&=\frac{|b+a''\cdot w_x''-w_{x,n}|}
        {\sqrt{1+|w_x''|^2}}\\
&=\frac{u_1|a\cdot x'+b-x_n|}
        {|u_1-x_1|\sqrt{1+|w_x''|^2}}\\
&\le \frac32|a\cdot x'+b-x_n|
\le  3nk\delta = \kappa.
\end{aligned}
\] 

\proofstep*{\textbf{Step 5:} Separation and the desired bound}

We now have both Frostman bounds and the required
incidence between points and hyperplanes. To finish, now express the point family $Z_2 \in H_{y_0}$ in the same $n-1$ coordinates:
\[
P_0=\{(a'',b):p=(a_1,a'',b)\in Z_2\}.
\]

The $\delta$-separated set $P_0$ still carries the required $(t-1)$-Frostman measure,
and by the preceding argument, each point of positive mass has an $s$-Frostman family
of hyperplanes meeting its $\kappa \delta$-ball. 

Let
\[
\mathcal T_0
=\bigcup_{p\in P_0}\mathcal T_p
\subset\{T_{\Pi_{y_0}(x)}:x\in E\}.
\]
Similarly, transfer weights to make $\mathcal T_0$ $\delta$-separated and keep the same notation.
By a finite partition and pigeonholing, we choose $\kappa\delta$-separated subsets
$
P\subset P_0,
\mathcal T\subset\mathcal T_0,
$
keeping a fixed positive proportion of the point mass and, at every remaining point, a positive proportion of its hyperplane probability. Restricting and normalizing preserves both Frostman bounds up to fixed factors related to $\kappa$.

Consequently, $(P,\mathcal T)$ is the required dual weighted
$(s,t-1)$ configuration in $\mathbb R^{n-1}$ at scale $\kappa\delta$. Moreover,
\[
\#\mathcal T\le\#\mathcal T_0\le\#E.
\]
Then we also conclude that 
\[
\mathcal N_{\kappa\delta}(\mathcal T)
\lesssim \mathcal N_\delta(E).
\]
Hence a bound for the covering number of the lower-dimensional hyperplane
family gives the same bound
for the original set $E$, up to a fixed factor.

\section{The terminal incidence theorem}

In this section, we restate the higher-dimensional hyperplane incidence
theorem of Cs\"ornyei and Stull and adapt it to our setting.
We use the full discretized statement of \cite[Theorem~1.1]{CS}
in $\mathbb R^n$, with the dyadic notation and definitions from
Section~\ref{sec:definitions}.

\begin{externaltheorem}[Cs\"ornyei--Stull, {\cite[Theorem~1.1]{CS}}]
Fix $n\ge2$, $\eps>0$, $n-2<s\le n-1$, and $0<t\le n$.
There exist $\eta>0$ and $\delta_0>0$ such that the following holds
for every dyadic $0<\delta<\delta_0$ and every $M\in\mathbb N$.
Let $\mathcal P\subseteq\mathcal D_\delta([0,1]^n)$ be a
$(\delta,t,\delta^{-\eta})$-set. For each $Q\in\mathcal P$, suppose
$H(Q)$ is a $(\delta,s,\delta^{-\eta})$-set of $\delta$-hyperplanes,
each intersecting $Q$, with $\#H(Q)\asymp M$.
Then $\mathcal H=\bigcup_{Q\in\mathcal P}H(Q)$ satisfies
\begin{equation}\label{eq:CS}
\#\mathcal H\gtrsim_\eps
M\delta^{-\min\{t,(s+t-n+2)/2,1\}+\eps}.
\end{equation}
\end{externaltheorem}

We now adapt this theorem to weighted configurations.

\begin{proposition}[Weighted terminal estimate]\label{prop:terminal}
Fix $n\ge2$ and bounded regions in $\R^n$ for the point set $E$ and
the hyperplane parameter set $Z$, independently of $\delta$.
For $n-2<s\le n-1$, $0<t\le n$ and every $\eps>0$, there exist
$\eta>0$, $c_\eps>0$ and $\delta_0>0$, depending only on
$n,s,t,\eps$ and these fixed regions, such that any weighted $(s,t)$
configuration whose points and hyperplane parameters lie in these regions,
with $K_s,K_t\le\delta^{-\eta}$ and $0<\delta<\delta_0$, satisfies
\begin{equation}\label{eq:weightedterminal}
\N_\delta(E)\ge c_\eps
\delta^{-\min\{s+t,(3s+t-n+2)/2,s+1\}+\eps}.
\end{equation}
Equivalently, for a dual weighted $(s,t)$ configuration the same
bound holds for $\N_\delta(\mathcal T)$, by
Definition~\ref{def:dualconfiguration}.
\end{proposition}

We first prove an extraction lemma needed for
Proposition~\ref{prop:terminal}.

\begin{lemma}[Extraction from a weighted measure]\label{lem:extract}
A Frostman $(\delta,\alpha,K)$ measure in a fixed bounded subset
of $\R^n$ gives a finite collection $\mathcal Q$ of occupied
$\delta$-cubes whose normalized Lebesgue measure on their union
is Frostman $(\delta,\alpha,C_2K\log(1/\delta))$, with a constant $C_2$ only depending on $n$, and
\[
\#\mathcal Q\ge
\bigl(C_2K\log(1/\delta)\delta^\alpha\bigr)^{-1}.
\]
\end{lemma}

\begin{proof}
There are at most $C_n\delta^{-n}$ occupied $\delta$-cubes $q$ meeting
the support of the given measure $\mu$. Discard those with
mass less than $\delta^{n+1}/(4C_n)$, losing at most $\delta/4$.
The remaining masses lie between $\delta^{n+1}/(4C_n)$ and $1$,
so there are only $O_n(\log(1/\delta))$ possible dyadic levels for pigeonholing.
By pigeonholing, there is a collection $\mathcal P'$ whose cube
masses lie in $[\lambda,2\lambda)$ for some $\lambda>0$, and whose
total mass is $\ge C_2 1/\log(1/\delta)$.
For $\delta\le r\le1$, counting the cubes meeting $B(x,r)$ gives
\[
\begin{aligned}
\#\{q\in\mathcal P':q \subset B(x,r)\}
&\le \frac{\mu(B(x,r))}{\lambda}\\
&\le \frac{Kr^\alpha}{\lambda}\\
&\le C_2 Kr^\alpha\log(1/\delta)
\frac{\mu\bigl(\bigcup_{q\in\mathcal P'}q\bigr)}{\lambda}\\
&\le C_2 Kr^\alpha\log(1/\delta)\,\#\mathcal P',
\end{aligned}
\]
which is the required Frostman bound for the normalised measure on the family of cubes $\mathcal P'$.
Taking $r=\delta$ and recomputing the above inequality gives
\[
\#\mathcal P'\ge (C_2K\log(1/\delta)\delta^\alpha)^{-1}.
\qedhere
\]
\end{proof}

\begin{proof}[Proof of Proposition~\ref{prop:terminal}]
After a fixed coordinate normalization, cover the hyperplane parameters by dyadic $\delta$-cubes and give each cube the sum of the original
weights. Applying Lemma~\ref{lem:extract} to this family measure gives
an unweighted Frostman
$
(\delta,t,C_2K_t\log(1/\delta))
$
set of dyadic $\delta$-hyperplanes. 

Choose one original hyperplane
$H_z$ in this family.
For each such $H_z$, apply the same extraction to $\mu_z$.
This gives an unweighted Frostman
$
(\delta,s,C_2K_s\log(1/\delta))
$
family $\mathcal P_z$ of $\delta$-cubes intersecting $H_z$, with
$
\#\mathcal P_z\ge
\bigl(C_2K_s\log(1/\delta)\delta^s\bigr)^{-1}.
$

We pigeonhole again because the Cs\"ornyei--Stull theorem requires the
number of cubes intersecting each hyperplane to be comparable to the same $M$.
By pigeonholing, keep a subfamily for which the values of $\# \mathcal P_z$ are comparable. Thus, for the remaining $H_z$,
\[
\#\mathcal P_z\asymp M,
\qquad
M\ge\bigl(C_2K_s\log(1/\delta)\delta^s\bigr)^{-1}.
\]
This costs one additional factor of $\log(1/\delta)$ in the Frostman
constant of the hyperplane family, which becomes $C_2K_t\log^2(1/\delta)$.
The Frostman constant of each $\mathcal P_z$ remains $C_2K_s\log(1/\delta)$.

Applying the Cs\"ornyei--Stull theorem after point--hyperplane duality,
with $\eta$ and then $\delta$ sufficiently small, gives
\[
\begin{aligned}
\N_\delta(E)
&\gtrsim_\eps
M\delta^{-\min\{t,(s+t-n+2)/2,1\}+\eps/2}\\
&\gtrsim_\eps
\frac{\delta^{-\min\{s+t,(3s+t-n+2)/2,s+1\}+\eps/2}}
{K_s\log(1/\delta)}.
\end{aligned}
\]
Since $K_s\le\delta^{-\eta}$, taking $\eta<\eps/4$ and $\delta$
sufficiently small gives the correct bound in the proposition.
\end{proof}

\section{Iteration and the Furstenberg bound}

For the iteration argument, write
\begin{equation}\label{eq:parameters}
s=d+\sigma,\qquad d=\lceil s\rceil-1,\qquad 0<\sigma\le1,
\qquad j=n-d-2.
\end{equation}

\begin{theorem}\label{thm:finite}
Fix $n\ge2$ and bounded regions in $\R^n$ for the point set $E$ and
the hyperplane parameter set $Z$, independently of $\delta$.
For $0<s\le n-1$, $0<t\le n$ and every $\eps>0$, there exist
$\eta>0$, $c_\eps>0$ and $\delta_0>0$, depending only on
$n,s,t,\eps$ and these bounded regions, such that any weighted $(s,t)$
configuration,
with $K_s,K_t\le\delta^{-\eta}$ and $0<\delta<\delta_0$, satisfies
\begin{equation}\label{eq:finite}
\N_\delta(E)\ge c_\eps\delta^{-\F(s,t;n,n-1)+\eps}.
\end{equation}
\end{theorem}
\begin{proof}
If $t\le j$, we can construct a
$t$-dimensional family of hyperplanes sharing a common
$\lceil s\rceil$-dimensional affine subspace, with the same
$s$-dimensional set of points lying on every hyperplane. This is
possible because the family of hyperplanes containing such a
subspace has dimension
\[
n-\lceil s\rceil-1=j.
\]
For any weighted configuration, one Frostman $(\delta,s,K_s)$
probability already gives
$
\N_\delta(E)\ge cK_s^{-1}\delta^{-s}.
$
Since $\F(s,t;n,n-1)=s$ in this range and $K_s\le\delta^{-\eta}$,
choosing $\eta<\eps$ proves \eqref{eq:finite}.

If $j=0$, apply Proposition~\ref{prop:terminal} directly.
Now suppose $j\ge1$ and $t>j$. Applying Theorem~\ref{thm:reduction}
gives a dual weighted $(s,t-1)$ configuration in $\R^{n-1}$, corresponding to the index
$\F(s,t-1;n-1,n-2)$. Dualizing back gives a weighted $(s,t-1)$
configuration in $\R^{n-1}$, to which we apply
Theorem~\ref{thm:reduction} again.
Repeating this process gives exactly $j$ reductions. 
At the step from ambient dimension $n-i$ to $n-i-1$, where
$0\le i<j$, we have
\[
s=d+\sigma\le d+1\le n-i-2,
\qquad t-i>j-i\ge1.
\]
Thus every step is admissible, including integer endpoints. The terminal
dimension $D$ and hyperplane family exponent $\tau$ are
\[
D=d+2 = \lceil s \rceil +1 , \qquad \tau =  t-j\in(0,D].
\]
In that dimension, notice that the exponent $s$ lies in the range where we are allowed to apply Proposition 4.1: 
$D-2<s\le D-1$. The accumulated Frostman constant losses are bounded by
$CK_s\log^{2j}(1/\delta)$ and $CK_t\log^{j}(1/\delta)$. Let \(\kappa_i\) be the scale multiplier at the \(i\)-th step. Each $\kappa_i$ is independent of $\delta$. The final configuration has scale
\[
\Delta=\left(\prod_{i=1}^{j}\kappa_i\right)\delta.
\]
Let $\eta_*>0$ be the exponent given by Proposition~\ref{prop:terminal}
in dimension $d+2$, with parameters $s,\tau$ and error $\eps$. Choose
$
{\eta=\min\left\{\frac{\eta_*}{2},\frac{\eps}{2}\right\}.}
$
Since \(\Delta\asymp\delta\) and \(\eta\le\eta_*/2\), for sufficiently small \(\delta\),
\[
K_s^{(j)},K_t^{(j)}
\le C\delta^{-\eta}\log^{2j}(1/\delta)
\le\Delta^{-\eta_*}.
\]
Thus Proposition~\ref{prop:terminal} applies to the final configuration at scale $\Delta$
and gives the required lower bound.
The restriction $\eta\le\eps/2$ ensures that in the end the loss is bounded by $\eps$. 
The exponent of the lower bound is
\begin{align}
\F(s,\tau;d+2,d+1)
&=\min\{s+\tau,(3s+\tau-d)/2,s+1\}\notag\\
&=s+\min\{\tau,(\sigma+\tau)/2,1\}\label{eq:matching}\\
&=\F(s,t;n,n-1).\notag
\end{align}
Combining the covering-number bounds from the reductions with Proposition 4.1 gives
\[
\begin{aligned}
\mathcal N_\delta(E)
&\gtrsim_\varepsilon
\Delta^{-\mathbf F(s,\tau;d+2,d+1)+\varepsilon}\\
&\asymp_n
\delta^{-\mathbf F(s,t;n,n-1)+\varepsilon}.
\end{aligned}
\]
This proves the required bound at the original scale \(\delta\).
\end{proof}

In particular, one admissible step has the numerical identity
\begin{equation}\label{eq:oneidentity}
\F(s,t;n,n-1)=\F(s,t-1;n-1,n-2)
\quad(0<s\le n-2,\ 1<t\le n).
\end{equation}

\section{The passage to Hausdorff dimension}

We prove Theorem~\ref{thm:hausdorff} using the covering argument of
\cite[Lemma~3.3]{HSY}, adapted to Theorem~\ref{thm:finite}.

\begin{proof}
Suppose, for a contradiction, that \(E\subset\mathbb R^n\) is Borel and \(\mathcal H\) is a family of hyperplanes satisfying
\[
\dim_H\mathcal H\ge t,
\qquad
\dim_H(E\cap H)\ge s
\quad(H\in\mathcal H),
\]
but
\[
\dim_HE<\mathbf F(s,t;n,n-1).
\]
Choose \(0<s_0<s\), \(0<t_0<t\), \(q>0\), and \(\varepsilon>0\), with \(\lceil s_0\rceil=\lceil s\rceil\), such that
\[
\dim_HE<q<\mathbf F(s_0,t_0;n,n-1)-\varepsilon.
\]
For integers \(R\ge1\), define
\[
\mathcal H_R
=
\{H\in\mathcal H:
\dim_H(E\cap H\cap[-R,R]^n)>s_0\}.
\]
Since \(\mathcal H=\bigcup_{R\ge1}\mathcal H_R\), countable stability gives an \(R\) with \(\dim_H\mathcal H_R>t_0\).
Replace $E$ by $E\cap [-R,R]^n$.
For the family $\mathcal H_R$, repeat the above procedure in parameter space. By countable stability, there exists an \(m\) such that
\[
\dim_H\{H_z\in\mathcal H_R:z\in[-m,m]^n\}>t_0.
\] After permuting coordinates, its hyperplanes have the form
\[
H_z=\{x:x_n=a\cdot x'+b\},
\qquad |a|\le mn,\qquad |b|\le m,
\]
where the bounds on both \(a \) and \(b\) follow because the parameters of the hyperplanes are bounded in \([-m,m]^n\).
Replace \(E\) by \(E\cap[-R,R]^n\) and \(\mathcal H\) by this subfamily, keeping the same notation. Both the points and hyperplane parameters now lie in fixed bounded regions. A lower bound for the dimension of this restricted set also gives the same lower bound for the original set.

We now have bounded points and parameters, with
\[
\mathcal H^{t_0}_\infty(\mathcal H)>0,
\qquad
\mathcal H^{s_0}_\infty(E\cap H)>0
\quad(H\in\mathcal H).
\]
Next, the subfamilies
\[
\mathcal H^{(M)}
=
\left\{
H\in\mathcal H:
\mathcal H^{s_0}_\infty(E\cap H)\ge\frac1M
\right\},
\qquad M\in\mathbb N,
\]
cover \(\mathcal H\). Countable subadditivity gives an \(M\) for which
\(\mathcal H^{t_0}_\infty(\mathcal H^{(M)})>0\). Put
\[
C_s=M,
\qquad
C_t=\frac1{\mathcal H^{t_0}_\infty(\mathcal H^{(M)})},
\]
and restrict to this subfamily. Thus
\[
\mathcal H^{t_0}_\infty(\mathcal H)=\frac1{C_t},
\qquad
\mathcal H^{s_0}_\infty(E\cap H)\ge\frac1{C_s}.
\]
We use the following finite-scale consequence of Frostman’s lemma: every bounded set \(A\) with positive \(\alpha\)-dimensional content has a probability on a finite \(\delta\)-separated subset satisfying
\[
\mu(B(x,r))
\lesssim
\frac{r^\alpha}{\mathcal H^\alpha_\infty(A)}
\qquad(\delta\le r\le1).
\]

Let \(\eta,\delta_0,c_\varepsilon\) be supplied by Theorem 5.1, decreasing \(c_\varepsilon\) by the fixed covering factor below if necessary. Choose \(k_0\) sufficiently large that
\[
\sum_{k\ge k_0}k^{-2}<1,
\qquad
C\max\{C_s,C_t\}k^2\le2^{k\eta}
\quad(k\ge k_0),
\]
and
\[
2^{-k_0}<\delta_0,
\qquad
c_\varepsilon
2^{k_0(\mathbf F(s_0,t_0;n,n-1)-q-\varepsilon)}>1.
\]
Since \(\dim_HE<q\), choose a cover of \(E\) by balls of dyadic radii smaller than \(2^{-k_0}\), with
$
\sum_i r_i^q<1.
$
Let \(\delta_k=2^{-k}\), let \(\mathcal Q_k\) be the balls of radius \(\delta_k\), and set
$
E_k=\bigcup_{Q\in\mathcal Q_k}Q.
$
Then
\[
E\subset\bigcup_{k\ge k_0}E_k,
\qquad
\sum_{k\ge k_0}\#\mathcal Q_k\,\delta_k^q<1.
\]
Define
\[
\mathcal H_k
=
\left\{
H\in\mathcal H:
\mathcal H^{s_0}_\infty(E\cap H\cap E_k)
\ge\frac1{C_sk^2}
\right\}.
\]
These families cover \(\mathcal H\). Indeed, a hyperplane belonging to none of them would satisfy
\[
\frac1{C_s}
\le
\mathcal H^{s_0}_\infty(E\cap H)
\le
\sum_{k\ge k_0}
\mathcal H^{s_0}_\infty(E\cap H\cap E_k)
\le
\frac1{C_s}\sum_{k\ge k_0}k^{-2}
<
\frac1{C_s}.
\]
Similarly,
\[
\frac1{C_t}
=
\mathcal H^{t_0}_\infty(\mathcal H)
\le
\sum_{k\ge k_0}\mathcal H^{t_0}_\infty(\mathcal H_k),
\]
so some \(k_1\ge k_0\) satisfies
\[
\mathcal H^{t_0}_\infty(\mathcal H_{k_1})
\ge\frac1{C_tk_1^2}.
\]
Put \(\delta=\delta_{k_1}\). The finite-scale Frostman lemma gives a finite \(\delta\)-separated family
\(\mathcal H'_{k_1}\subset\mathcal H_{k_1}\) carrying a probability \(\nu\), and, for each \(H\in\mathcal H'_{k_1}\), a finitely supported probability \(\mu_H\) on \(E\cap H\cap E_{k_1}\), with
\[
\nu(B(z,r))\lesssim C_tk_1^2r^{t_0},
\qquad
\mu_H(B(x,r))\lesssim C_sk_1^2r^{s_0}
\quad(\delta\le r\le1).
\]
Only finitely many measures are being chosen.
Choose a maximal \(\delta\)-separated subset \(E'_{k_1}\subset E\cap E_{k_1}\), and transfer each \(\mu_H\) to nearest representatives. The resulting probability is supported within distance \(\delta\) of \(H\), and
\[
\mu_{H,\delta}(B(x,r))
\le\mu_H(B(x,r+\delta))
\lesssim C_sk_1^2r^{s_0}.
\]
We therefore have a weighted \((s_0,t_0)\) configuration with
\[
K_s\lesssim C_sk_1^2,
\qquad K_t\lesssim C_tk_1^2,
\qquad
\mathcal N_\delta(E'_{k_1})\lesssim\#\mathcal Q_{k_1}.
\]
By the choice of \(k_0\), Theorem 5.1 applies and gives
\[
\#\mathcal Q_{k_1}
\ge
c_\varepsilon
\delta_{k_1}^{-\mathbf F(s_0,t_0;n,n-1)+\varepsilon}.
\]
Consequently,
\[
\begin{aligned}
\sum_{k\ge k_0}\#\mathcal Q_k\,\delta_k^q
&\ge \#\mathcal Q_{k_1}\delta_{k_1}^q\\
&\ge
c_\varepsilon
2^{k_1(\mathbf F(s_0,t_0;n,n-1)-q-\varepsilon)}
>1,
\end{aligned}
\]
contradicting the choice of the cover. This proves Theorem~\ref{thm:hausdorff}.
\end{proof}

\section{Sharpness and the higher-dimensional examples}\label{sec:sharpness}

For every $0<s\le n-1$ and $0<t\le n$, we give compact sets
$E\subset\R^n$ and compact hyperplane families $\HH$ such that
\[
\dimH\HH=t,\qquad \dimH(E\cap H)=s\quad(H\in\HH),
\qquad \dimH E=\F(s,t;n,n-1).
\]
Thus the lower bound in Theorem~\ref{thm:hausdorff} is sharp.
We use the notation from Section~5:
\[
s=d+\sigma,\qquad 0<\sigma\le1,\qquad j=n-d-2.
\]
Let $C_\sigma \subset[0,1]$ be a compact AD-regular set of dimension
$\sigma$, using an interval when $\sigma=1$.

The sharp examples in $2D$ are already known; see
\cite[Section~1]{RW} and \cite{Wolff}. For every $0<\sigma\le1$
and $0<\tau\le2$, we can find a compact set $E_0\subset\R^2$
and a compact family $\mathcal L$ of lines
$\ell_{a,b}=\{(x,y):y=ax+b\}$ such that
\[
\dimH\mathcal L=\tau,\qquad
\dimH(E_0\cap\ell)=\sigma\quad(\ell\in\mathcal L),
\]
and
\[
\dimH E_0=\min\{\sigma+\tau,(3\sigma+\tau)/2,\sigma+1\}.
\]
We use these examples to construct the examples in $\R^n$ below.

\subsection{The examples in $\R^n$}

\paragraph{If $0<t\le j$.} Sharpness follows from the construction in Héra’s Proposition 1.7, specialized to hyperplanes and restricted to a $(t)$-dimensional subfamily.
Place $E=C_\sigma\times[0,1]^d$ in a fixed $(d+1)$-dimensional
linear subspace of $\R^n$. The family of all hyperplanes containing
this subspace has dimension $n-d-2=j$.
Choose a compact $t$-dimensional subfamily $\HH$ inside the set.
Every $H\in\HH$ contains $E$, and therefore
\[
\dimH(E\cap H)=\dimH E=d+\sigma=s=\F(s,t;n,n-1).
\]

\paragraph{If $j<t<j+2-\sigma$.} The example extends the planar $2D$ case and adds dimensions up.
Put $\tau=t-j$ and take the planar example $(E_0,\mathcal L)$ above, where $\dimH\mathcal L=\tau$.
Write coordinates in $\R^n$ as
$(x,y,v,w)\in\R\times \R\times\R^d\times\R^j$, and set
\[
E=E_0\times[0,1]^d\times\{0\}^j,\qquad
H_{a,b,c}=\{(x,y,v,w):y=ax+b+c\cdot w\},
\]
where $\ell_{a,b} = \{(x,y) \in \mathbb R^2: y = ax +b\}\in\mathcal L$ and $c\in[0,1]^j$.
The hyperplane parameter family has dimension $\tau+j=t$.
Also,
\[
E\cap H_{a,b,c}
=(E_0\cap\ell_{a,b})\times[0,1]^d\times\{0\}^j,
\]
so every intersection has dimension $\sigma+d=s$. Finally,
\[
\begin{aligned}
\dimH E
&=d+\min\{\sigma+\tau,(3\sigma+\tau)/2\}\\
&=\min\{s+t-j,(3s+t-n+2)/2\}
=\F(s,t;n,n-1).
\end{aligned}
\]

\paragraph{If $j+2-\sigma\le t\le n$.} We can simplify the example above since the resulting Furstenberg dimension is always $s+1$. The overall construction is very similar to the sharp example in the $2$-dimensional case when $s+t>2$.
In the same coordinates, take
\[
E=C_\sigma\times[0,1]\times[0,1]^d\times\{0\}^j.
\]
Choose a compact \(t\)-dimensional set \(A\) of parameters \((a,b,\xi,c) \in \R \times \R \times \R^{d}\times \R^j\) in $\R^n$ satisfying
\[
|a|,|\xi_i|\le\frac1{8(d+1)},\qquad
\frac14\le b\le\frac34,\qquad c\in[0,1]^j,
\]
and define
\[
\mathcal H=
\left\{
H_{a,b,\xi,c}=\{y=ax+b+\xi\cdot v+c\cdot w\}:
(a,b,\xi,c)\in A
\right\}.
\]
By duality, 
\(\dim_H\mathcal H=t\).
For \(x\in C_\sigma\), \(v\in[0,1]^d\), and the last coordinate \(w=0\), we have
\[
\frac18\le ax+b+\xi\cdot v\le\frac78.
\]
Thus each \(E\cap H\) is a graph over \(C_\sigma\times[0,1]^d\) and has dimension \(s\). Since \(t-j\ge2-\sigma\),
\[
\dim_HE=d+\sigma+1=s+1
=\mathbf F(s,t;n,n-1).
\]

\subsection{The plane Furstenberg problem in $\R^3$}\label{sec:three-dimensional}

In this final subsection, we present the explicit formula for the index when the hyperplanes are two-dimensional planes and give some intuition for the bound, with two cases
$0<s\le1$ and $1<s\le2$. The sharp bound is
$\dimH E\ge\F(s,t;3,2)$, and the exact parameter ranges are as follows. In higher dimensions, the formulas get more complicated, but the same induction process applies.
\paragraph{If $0<s\le1$.} The condition \(s\le1\) allows us to apply radial reduction. For \(t>1\), one radial reduction gives a weighted \((s,t-1)\) configuration in \(\mathbb R^2\), to which the planar Ren–Wang estimate applies. The last three expressions below follow by replacing \(t\) with \(t-1\) in the planar bound:
\[
\min\left\{s+t-1,\frac{3s+t-1}{2},s+1\right\}.
\]
For $t\le1$, all planes can contain the same line, on which we place an $s$-dimensional set. Thus
\[
\F(s,t;3,2)=
\begin{cases}
s,&0<t\le1,\\[2pt]
s+t-1,&1<t\le1+s,\\[2pt]
(3s+t-1)/2,&1+s<t\le3-s,\\[2pt]
s+1,&3-s<t\le3.
\end{cases}
\]

\paragraph{If $1<s\le2$.} 
The Cs\"ornyei--Stull theorem applies directly, giving
\[
\F(s,t;3,2)=
\begin{cases}
s+t,&0<t\le s-1,\\[2pt]
(3s+t-1)/2,&s-1<t\le3-s,\\[2pt]
s+1,&3-s<t\le3.
\end{cases}
\]
For \(0<t\le2\), this bound can be understood through a planar product construction. Start with an \((s-1,t)\)-Furstenberg example \(E_0\) in the \(xy\)-plane. Extend each line \(\ell\) in its family in the \(z\)-direction to obtain the plane \(\ell\times\mathbb R\), and set
$
E=E_0\times[0,1].
$
Each resulting plane is perpendicular to the \(xy\)-plane, and its intersection with \(E\) has dimension \(s\). Taking the product with an interval increases the dimension of \(E_0\) by one. Thus, replacing \(s\) by \(s-1\) in the planar index and adding one gives
\[
\begin{aligned}
1+\min\left\{s-1+t,\frac{3(s-1)+t}{2},s\right\}
&=\min\left\{s+t,\frac{3s+t-1}{2},s+1\right\}.
\end{aligned}
\]
For \(2<t\le3\), the index is \(s+1\), and the preceding \(s+1\)-dimensional example provides the required hyperplane family.
In particular, in both cases,
$s+t\ge3$ gives $\F(s,t;3,2)=s+1$.

In three dimensions, the author also obtained the same bound in
the range $s+t>3$ by adapting the high--low argument of
Fu and Ren~\cite{FR}. The author separately proved the plane case
$1<s<2$ by combining the Ren--Wang argument~\cite{RW} with the
$t$-AD-regular case. These proofs follow the central ideas of the
cited papers. In view of the cleaner result of Cs\"ornyei--Stull,
these alternative proofs are not reproduced here.

\end{document}